\documentclass[11pt, notitlepage]{article}
\usepackage{amssymb,amsmath,comment}
\usepackage{amsthm}
\catcode`\@=11 \@addtoreset{equation}{section}
\def\thesection{\arabic{section}}

\def\theequation{\thesection.\arabic{equation}}
\catcode`\@=12
\usepackage{colortbl}
\usepackage{hyperref}
\usepackage[mathscr]{eucal}
\usepackage{epsf}
\usepackage{esint}
\usepackage{a4wide}
\def\R{\mathbb{R}}

\DeclareMathOperator*{\esssup}{ess\,sup}
\DeclareMathOperator*{\essinf}{ess\,inf}

\DeclareMathOperator*{\supp}{supp}

\newcommand{\Om} {\Omega}

\newcommand{\noi} {\noindent}

\newcommand{\Tail} {\mathrm{Tail}}
\newcommand{\loc} {\mathrm{loc}}

\usepackage[all]{xy}
\catcode`\@=11
\def\theequation{\@arabic{\c@section}.\@arabic{\c@equation}}
\catcode`\@=12

\newtheorem{Theorem}{Theorem}[section]
\newtheorem{Lemma}[Theorem]{Lemma}

\newtheorem{Remark}[Theorem]{Remark}
\newtheorem{Definition}[Theorem]{Definition}

\begin{document}

{\vspace{0.01in}}

\title{On the regularity theory for mixed anisotropic and nonlocal $p$-Laplace equations and its applications to singular problems}
\author{Prashanta Garain, Wontae Kim, Juha Kinnunen}

\date{}
\maketitle

\begin{abstract}
We establish existence results for a class of mixed anisotropic and nonlocal $p$-Laplace equation with singular nonlinearities. We consider both constant and variable singular exponents. Our argument is based on  an approximation method. To this end, we also discuss the necessary regularity properties of weak solutions of the associated non-singular problems. More precisely, we obtain local boundedness of subsolutions, Harnack inequality for solutions and weak Harnack inequality for supersolutions.
\end{abstract}

\maketitle

\noi {Keywords: Mixed anisotropic and nonlocal $p$-Laplace equation, existence, regularity, singular problem, variable exponent.}

\noi{\textit{2020 Mathematics Subject Classification: 35M10, 35R11, 35B65, 35J75, 35J92.}

\section{Introduction}
This article discusses regularity properties of weak solutions of a mixed local and nonlocal $p$-Laplace equation
\begin{equation}\label{genmaineqn}
-\operatorname{div}\mathcal{B}(x,u,\nabla u)+\mathcal{L}(u)=0\text{ in }\Omega
\end{equation}
in a bounded smooth domain $\Omega$ of $\mathbb{R}^N$. Here $\mathcal{B}(x,u,\zeta):\Om\times\mathbb{R}^{N+1}\to\mathbb{R}^N$ is a measurable function with respect to $x$ and continuous with respect to $(u,\zeta)$ such that 
\begin{itemize}
\item[(H1)] $\mathcal{B}(x,u,\zeta)\cdot\zeta\geq C_1|\zeta|^p$,\text{ and }$|\mathcal{B}(x,u,\zeta)|\leq C_2|\zeta|^{p-1}$
for almost every $x\in\Om$ and for every $\,(u,\zeta)\in\mathbb{R}^{N+1}$,
where $C_1,C_2$ are some positive constants and $1<p<\infty$.
\end{itemize}
The operator $\mathcal{L}$ is the nonlocal $p$-Laplace operator defined by
\begin{equation}\label{flap}
\mathcal{L}(u)(x)=\text{P.V.}\int_{\mathbb{R}^N}|u(x)-u(y)|^{p-2}(u(x)-u(y))K(x,y)\,dy,
\end{equation}
where P.V. denotes the principal value and $K$ is a symmetric kernel in $x$ and $y$ such that
\begin{equation}\label{flapker}
\frac{\Lambda^{-1}}{|x-y|^{N+ps}}\leq K(x,y)\leq \frac{\Lambda}{|x-y|^{N+ps}}
\end{equation}
for some constant $\Lambda\geq 1$ and $0<s<1$.

We apply the obtained regularity results to investigate existence results for the Dirichlet problems for a mixed local and nonlocal $p$-Laplace equation with singular nonlinearities (both constant and variable exponent)
\begin{equation}\label{sinmeqn1}
-\text{div}(\mathcal{B}(x,\nabla u))+\mathcal{L}(u)=f(x)u^{-\gamma(x)}\text{ in }\Omega,\quad u>0\text{ in }\Om,\quad u=0\text{ in }\mathbb{R}^N\setminus\Omega,
\end{equation}
where $\mathcal{B}(x,\nabla u)$ satisfy (H1) along with the strict monotonicity condition
\begin{itemize}
\item[(H2)] $(\mathcal{B}(x,\zeta_1)-\mathcal{B}(x,\zeta_2))(\zeta_1-\zeta_2)>0$ for almost every $x\in\Om,$ and for every $\zeta_1,\zeta_2\in\mathbb{R}^N$ with $\zeta_1\neq \zeta_2$.
\end{itemize}
There are several interesting equations of type \eqref{genmaineqn}. One of the leading examples is a mixed local and nonlocal $p$-Laplace equation
\begin{equation}\label{p1}
-a\Delta_p u+b(-\Delta_p)^s u=0
\end{equation}
for some positive constants $a,b$, which is obtained by choosing 
$$
\mathcal{B}(x,u,\zeta)=a|\zeta|^{p-2}\zeta,\quad\text{and}\quad K(x,y)=b|x-y|^{-N-ps}.
$$
The operators $\Delta_p$ and $(-\Delta_p)^s$ are the $p$-Laplace and the fractional $p$-Laplace operators respectively.
A more general mixed anisotropic and nonlocal $p$-Laplace equation
\begin{equation}\label{p11}
-aH_p u+b(-\Delta_p)^s u=0
\end{equation}
can be obtained from \eqref{genmaineqn} by choosing 
$$
\mathcal{B}(x,u,\zeta)=aH(\zeta)^{p-1}\nabla H(\zeta),\quad\text{and}\quad K(x,y)=b|x-y|^{-N-ps}.
$$
The operator
\begin{equation}\label{fo}
H_p u=\text{div}\big(H(\nabla u)^{p-1}\nabla H(\nabla u)\big),
\end{equation}
is referred to as the anisotropic $p$-Laplace operator, where $H:\R^N\to[0,\infty)$ is a Finsler-Minkowski norm,  that is, $H$ is a nonnegative convex function in $C^1(\R^N\setminus\{0\})$ such that $H(\zeta)=0$ if and only if $\zeta=0$, and $H$ is even and positively homogeneous of degree $1$, so that
$
H(t\zeta)=|t|H(\zeta)
$
\text{for every $\zeta\in\R^N$ and $t\in\mathbb{R}$}. Then from Xia \cite[Chapter $1$]{Xiathesis} we observe that $\mathcal{B}(x,u,\zeta)=H(\zeta)^{p-1}\nabla H(\zeta)$ satisfies (H1) for every $1<p<\infty$. Furthermore, from \cite[Lemma 2.5]{BGM}, it follows that $\mathcal{B}(x,u,\zeta)=H(\zeta)^{p-1}\nabla H(\zeta)$ satisfies (H2) for every $2\leq p<\infty$.
Various examples of Finsler-Minkowski norm $H$ can be found in, for example in Belloni-Ferone-Kawohl \cite{BFKzamp}, Xia \cite[p. 22--23]{Xiathesis} and the references therein. A typical example of $H$ includes the $l^q$-norm defined by
\begin{equation}\label{ex1}
H(\zeta)=\Big(\sum_{i=1}^{N}|\zeta_i|^q\Big)^\frac{1}{q},\quad q>1,
\end{equation}
where $\zeta=(\zeta_1,\zeta_2,\ldots,\zeta_N)$. When $H$ is the $l^q$-norm as in \eqref{ex1}, we have
\begin{equation}\label{exf}
H_p u=\sum_{i=1}^{N}\frac{\partial}{\partial x_i}\bigg(\Big(\sum_{k=1}^{N}\Big|\frac{\partial u}{\partial x_k}\Big|^{q}\Big)^\frac{p-q}{q}\Big|\frac{\partial u}{\partial x_i}\Big|^{q-2}\frac{\partial u}{\partial x_i}\bigg).
\end{equation}
For $q=2$ in \eqref{exf}, $H_p$ becomes the usual $p$-Laplace operator $\Delta_p$. Moreover, for $q=p$ in \eqref{exf}, the operator $H_p$ reduces to the pseudo $p$-Laplace operator, see Belloni-Kawohl \cite{BKesaim} and therefore, equation \eqref{genmaineqn} covers a wide range of mixed local and nonlocal problems and, in particular, extends the following mixed equation
\begin{equation}\label{mpf}
\sum_{i=1}^{N}\frac{\partial}{\partial x_i}\bigg(\Big(\sum_{k=1}^{N}\Big|\frac{\partial u}{\partial x_k}\Big|^{q}\Big)^\frac{p-q}{q}\Big|\frac{\partial u}{\partial x_i}\Big|^{q-2}\frac{\partial u}{\partial x_i}\bigg)+(-\Delta_p)^s u=0
\text{ in }\Om.
\end{equation}
When $2\leq p<\infty$, the singular equation in \eqref{sinmeqn1} extends the following mixed anisotropic and fractional $p$-Laplace equation 
\begin{equation}\label{sinmeqn2}
-aH_p u+b(-\Delta_p)^s u=f(x)u^{-\gamma(x)}\text{ in }\Omega,\quad u>0\text{ in }\Om,\quad u=0\text{ in }\mathbb{R}^N\setminus\Omega.
\end{equation}

The singular equation in \eqref{sinmeqn2} has been studied extensively when $a=1,\,b=0$ (local case) and $a=0,\,b=1$ (nonlocal case). 
In the local case, the singular $p$-Laplace equation
\begin{equation}\label{splap}
-\Delta_p u=u^{-\gamma},\quad u>0\text{ in }\Om,\quad u=0\text{ on }\partial\Om
\end{equation}
has been studied thoroughly over the last three decades. When $\gamma>0$ is a constant, for $p=2$, Crandall-Rabinowitz-Tartar \cite{CRT} proved the existence of a classical solution for any $\gamma>0$. Lazer-Mckenna \cite{LazMc} proved such a solution is weak if and only if $0<\gamma<3$. This restriction on $\gamma$ is removed by Boccardo-Orsina \cite{BocOr} to obtain weak solutions for any $\gamma>0$. 
Furthermore, equation in \eqref{splap} is studied to obtain existence and uniqueness of weak solutions for $1<p<\infty$ and $\gamma>0$ by Canino-Trombetta-Sciunzi \cite{Canino}. 

When $\gamma$ is a variable, for $p=2$, existence results for \eqref{splap} have been obtained by Carmona-Mart\'{\i}nez-Aparicio \cite{CMP} and further extended to the variants of $p$-Laplace equations by Alves-Santos-Siqueira \cite{Alves20}, Alves-Moussaoui \cite{Alves18}, Byun-Ko \cite{BKcvpde}, Papageorgiou-Scapellato \cite{PSzamp}, Chu-Gao-Gao \cite{CGG}, Zhang \cite{Zhang}.
Recently singular problems have also been studied in the anisotropic setting. For a constant singular exponent, we refer to Biset-Mebrate-Mohammed \cite{BMM20}, Farkas-Winkert \cite{PF20}, Farkas-Fiscella-Winkert \cite{PF21}, Garain \cite{Gtmna}, Miri \cite{Mirivar} and for a variable singular exponent, see Bal-Garain-Mukherjee \cite{BGM, Ganisovar} and the references therein.

In the nonlocal case, the Dirichlet problems for a singular fractional $p$-Laplace equation
\begin{equation}\label{splap1}
(-\Delta_p)^s u=u^{-\gamma},\quad u>0\text{ in }\Om,\quad u=0\text{ in }\mathbb{R}^N\setminus\Om
\end{equation}
is also studied widely. For a  constant singular exponent $\gamma$, we refer to Fang \cite{Fang} in the semilinear case $p=2$; Canino-Montoro-Sciunzi-Squassina \cite{Caninononloc} in the quasilinear case. The perturbed  singular case is investigated by Barrios-De Bonis-Medina-Peral \cite{BDMP}, Adimurthi-Giacomoni-Santra \cite{Adi} for $p=2$; Mukherjee-Sreenadh \cite{MS} in the quasilinear case and the references therein. For a variable exponent $\gamma$, see Garain-Mukherjee \cite{GMnonloc} in the quasilinear case $1<p<\infty$. We also refer to the nice survey by Ghergu-R\u{a}dulescu \cite{GRbook} on singular elliptic problems.

Mixed singular problems are much less understood. Recently, for a constant singular exponent $\gamma$, the semilinear mixed local and nonlocal Laplace equation
\begin{equation}\label{mlap}
-\Delta u+(-\Delta)^s u={f(x)}{u^{-\gamma}},\quad u>0\text{ in }\Om,\quad u=0\text{ in }\mathbb{R}^N\setminus\Om
\end{equation}
is investigated by Arora-Radulescu \cite{AroRad} to prove existence results. We also refer to Garain \cite{Gjga} for semilinear mixed probelms with purturbed singularity. The quasilinear case $1<p<\infty$ for the equation
\begin{equation}\label{mlap1}
-\Delta_p u+(-\Delta_p)^s u=f(x){u^{-\gamma}},\quad u>0\text{ in }\Om,\quad u=0\text{ in }\mathbb{R}^N\setminus\Om
\end{equation}
is studied by Garain-Ukhlov \cite{GU} to obtain existence, uniqueness and regularity results for the equation \eqref{mlap1}. Very recently, for variable singular exponent $\gamma\in C^1(\overline{\Omega})$, Biroud \cite{Biroud} studied existence results for the mixed $p$-Laplace equation \eqref{mlap1} for $1<p<N$. In this article, we consider the mixed anistropic and nonlocal $p$-Laplace equation \eqref{genmaineqn} which extends the equation \eqref{mlap1}. Further, we only assume $\gamma\in C(\overline{\Omega})$ and discuss existence results in the whole range $1<p<\infty$. Moreover, we investigate the case of constant singular exponent $\gamma$ which improves the range of the source term $f$ to obtain existence results as stated in Section 4.

\section{Preliminaries}
The following notation will be used throughout the paper.
\begin{itemize}
\item For $k\in\mathbb{R}$, we denote by $k^+=\max\{k,0\}$, $k^-=\max\{-k,0\}$ and $k_-=\min\{k,0\}$.
\item $\langle,\rangle$ denotes the standard inner product in $\mathbb{R}^N$.
\item For $r>1$, we denote by $r'=\frac{r}{r-1}$ to mean the conjugate exponent of $r$.
\item For $1<p<N$, we denote by $p^*=\frac{Np}{N-p}$ to mean the critical Sobolev exponent.
\item For a function $F$ defined over a set $S$ and some constants $c$ and $d$, by $c\leq F\leq d$ in $S$, we mean $c\leq F\leq d$ almost everywhere in $S$.
\item We denote an open ball with center $x_0\in\R^n$ and radius $r>0$ by $B_r(x_0)$.
\item $\Omega$ is a bounded smooth domain $\mathbb{R}^N$ and 
\[
\Omega_{\delta}=\{x\in\Om:\text{dist}(x,\partial\Om)<\delta\},\quad\delta>0.
\]
\item $C$ will denote a constant which may vary from line to line or even in the same line. If $C$ depends on the parameters $r_1,r_2,\ldots$,  we write $C=C(r_1,r_2,\ldots,)$.
\end{itemize}

In this section, we present some known results for the fractional Sobolev spaces, see \cite{Hitchhiker'sguide} for more details.
Let 
$E\subset \mathbb{R}^N$
be a measurable set and $|E|$ denote its Lebesgue measure. Recall that the Lebesgue space 
$L^{p}(E),1\leq p<\infty$, 
is defined as the space of $p$-integrable functions $u:E\to\mathbb{R}$ with the norm
$$ \|u\|_{L^p(E)}=
\left(\int_{E}|u(x)|^p~dx\right)^{\frac1p}.
$$
By 
$L^p_{\mathrm{loc}}(E)$ 
we denote the space of locally $p$-integrable functions, means that,
$u\in L^p_{\mathrm{loc}}(E)$ 
if and only if  
$u\in L^p(F)$ 
for every compact subset 
$F\subset E$. In the case $0<p<1$, we denote by $L^{p}(E)$ a set of measurable functions such that 
$\int_{E}|u(x)|^p~dx<\infty$.
Here and in the rest of the paper, let $\Om\subset\mathbb{R}^N$ with $N\geq 2$ be a bounded smooth domain. 
The Sobolev space $W^{1,p}(\Omega)$, $1<p<\infty$, is defined 
as the Banach space of locally integrable weakly differentiable functions
$u:\Omega\to\mathbb{R}$ equipped with the norm
\[
\|u\|_{W^{1,p}(\Omega)}=\| u\|_{L^p(\Omega)}+\|\nabla u\|_{L^p(\Omega)}.
\]
The space $W^{1,p}(\mathbb{R}^N)$ is defined analogously.
The Sobolev space with zero boundary values is defined as
\[
W^{1,p}_0(\Omega)=\{u\in W^{1,p}(\mathbb{R}^N):u=0\text{ in }\mathbb{R}^N\setminus\Om\}
\]
under the norm $\|u\|=\|\nabla u\|_{L^p(\Om)}$. 
This space is a separable and reflexive Banach space for any $1<p<\infty$, see \cite{BDVV2, Biagi1, SV}.  

The fractional Sobolev space $W^{s,p}(\Omega)$, $0<s<1<p<\infty$, is defined by
$$
W^{s,p}(\Omega)=\Big\{u\in L^p(\Omega):\frac{|u(x)-u(y)|}{|x-y|^{\frac{N}{p}+s}}\in L^p(\Omega\times \Omega)\Big\}
$$
and endowed with the norm
$$
\|u\|_{W^{s,p}(\Omega)}=\left(\int_{\Omega}|u(x)|^p\,dx+\int_{\Omega}\int_{\Omega}\frac{|u(x)-u(y)|^p}{|x-y|^{N+ps}}\,dx\,dy\right)^\frac{1}{p}.
$$
The space $W^{s,p}_{\mathrm{loc}}(\Omega)$ is defined analogously.
The next result asserts that the Sobolev space $W^{1,p}(\Om)$ is continuously embedded in the fractional Sobolev space, see \cite[Proposition $2.2$]{Hitchhiker'sguide}.
\begin{Lemma}\label{defineq}
There exists a constant $C=C(N,p,s)$ such that
$\|u\|_{W^{s,p}(\Om)}\leq C\|u\|_{W^{1,p}(\Om)}$
for every $u\in W^{1,p}(\Omega)$.
\end{Lemma}

Next, we have the following result from \cite[Lemma $2.1$]{Silva}.
\begin{Lemma}\label{locnon1}
There exists a constant $c=c(N,p,s,\Omega)$ such that
\begin{equation}\label{locnonsem}
\int_{\mathbb{R}^N}\int_{\mathbb{R}^N}\frac{|u(x)-u(y)|^p}{|x-y|^{N+ps}}\,dx\,dy\leq c\int_{\Omega}|\nabla u|^p\,dx
\end{equation}
for every $u\in W_0^{1,p}(\Omega)$.
\end{Lemma}

For the following Sobolev embedding, see, for example, \cite{Evans}.

\begin{Lemma}\label{emb}
The embedding operators
\[
W_0^{1,p}(\Omega)\hookrightarrow
\begin{cases}
L^t(\Om),&\text{ for }t\in[1,p^{*}],\text{ if }1<p<N,\\
L^t(\Om),&\text{ for }t\in[1,\infty),\text{ if }p=N,\\
L^{\infty}(\Om),&\text{ if }p>N,
\end{cases}
\]
are continuous. Also, the above embeddings are compact, except for $t=p^*=\frac{Np}{N-p}$, if $1<p<N$.
\end{Lemma}

Next, we state the algebraic inequality from \cite[Lemma $2.1$]{Dama}.

\begin{Lemma}\label{alg}
Let $1<p<\infty$. Then for any $a,b\in\mathbb{R}^N$, there exists a constant $c=c(p)>0$
such that
\begin{equation}\label{algineq}
(|a|^{p-2}a-|b|^{p-2}b)(a-b) \rangle\geq
c\frac{|a-b|^2}{(|a|+|b|)^{2-p}}.
\end{equation}
\end{Lemma}

The following result follows from \cite[Lemma A.2]{BrPr}.

\begin{Lemma}\label{BrPrapp}
Let $1<p<\infty$ and $g:\mathbb{R}\to\mathbb{R}$ be an increasing function. Then for every $a,b\in\mathbb{R}$, we have
\begin{equation}\label{BrPrine}
|a-b|^{p-2}(a-b)\big(g(a)-g(b)\big)\geq|G(a)-G(b)|^p,
\end{equation}
where
$$
G(t)=\int_{0}^{t}g'(\tau)^\frac{1}{p}\,d\tau,\quad t\in\mathbb{R}.
$$
\end{Lemma}

The following result from \cite[Theorem $9.14$]{var} will be useful to prove Lemma \ref{exisapprox} below.

\begin{Theorem}\label{MB}
Let $V$ be a real separable reflexive Banach space and $V'$ be the dual of $V$. Suppose that $T:V\to V'$ is a coercive and demicontinuous monotone operator. Then $T$ is surjective, i.e., given any $f\in V'$, there exists $u\in V$ such that $T(u)=f$. If $T$ is strictly monotone, then $T$ is also injective.  
\end{Theorem}

Next we will define the notions of a weak solution of \eqref{genmaineqn} and \eqref{sinmeqn1}.

\begin{Definition}\label{gensubsupsolution}
Let $1<p<\infty$.
A function $u\in L^{\infty}(\mathbb{R}^N)$ is a weak subsolution of \eqref{genmaineqn} if $u\in W_{\loc}^{1,p}(\Omega)$ and for every $\Omega'\Subset\Omega$ and every nonnegative test function $\phi\in W_{0}^{1,p}(\Omega')$, we have
\begin{equation}\label{weaksubsupsoln}
\begin{gathered}
\int_{\Omega'}\mathcal{B}(x,u,\nabla u)\cdot\nabla\phi(x)\,dx+\int_{\mathbb{R}^n}\int_{\mathbb{R}^n}\mathcal{A}(u(x,y)){(\phi(x)-\phi(y))}\,d\mu\leq0,
\end{gathered}
\end{equation}
where
\begin{equation}\label{not}
\mathcal{A}(u(x,y))=|u(x)-u(y)|^{p-2}(u(x)-u(y)),
\quad
d\mu=K(x,y)\,dx\,dy
\end{equation}
and $K$ is defined in \eqref{flapker}.
Analogously, a function $u$  is a weak supersolution of  \eqref{genmaineqn} if the integral in \eqref{weaksubsupsoln} is nonnegative for every nonnegative test function $\phi\in W_{0}^{1,p}(\Omega')$.
A function $u$ is a weak solution of \eqref{genmaineqn} if the equality holds in \eqref{weaksubsupsoln} for every $\phi\in W_{0}^{1,p}(\Omega')$ without a sign restriction.
\end{Definition}

\begin{Remark}\label{rkreg}The boundedness assumption together with Lemma \ref{defineq} and Lemma \ref{locnon1}, ensures that  Definition \ref{gensubsupsolution} is well stated and
 the tail defined in \eqref{loctail} is finite. 
Under the assumption that the tail in \eqref{loctail} is bounded, our main regularity results Theorem \ref{locbdd}, Theorem \ref{regthm2} and Theorem \ref{regthm3} hold true without the a priori boundedness assumption on the function. 
In this case, the local boundedness follows from Theorem \ref{locbdd}.
\end{Remark}

We define the notion of zero Dirichlet boundary condition as follows:

\begin{Definition}\label{bc}
We say that $u\leq 0$ on $\partial\Om$, if $u=0\text{ in }\mathbb{R}^N\setminus\Om$ and for every $\epsilon>0$, we have
$$
(u-\epsilon)^+\in W_0^{1,p}(\Om).
$$
We say that $u=0$ on $\partial\Om$, if $u$ is nonnegative and $u\leq 0$ on $\partial\Om$.
\end{Definition}

\begin{Remark}\label{bcrmk}
If $u\in W^{1,p}_{\mathrm{loc}}(\Om)$ is nonnegative in $\Om$, such that $u^\alpha\in W_0^{1,p}(\Om)$, for some $\alpha\geq 1$, then following the proof of
{\rm\cite[Page 5]{Canino}} or {\rm\cite[Theorem 2.11]{G}} we have $u=0$ on $\partial\Om$ according to the Definition \ref{bc}. 
\end{Remark}

Next, we define the notion of weak solution for  \eqref{sinmeqn1}.

\begin{Definition}\label{psintest}
We say that a function $u\in W_{\mathrm{loc}}^{1,p}(\Omega)\cap L^{p-1}(\Om)$ is a weak solution of \eqref{sinmeqn1}, if
$u>0$ in $\Om$, $u=0$ on $\partial\Om$ in the sense of Definition \ref{bc} and  ${f}{u^{-\gamma(x)}}\in L^1_{\mathrm{loc}}(\Om)$,
such that for every $\phi\in C_{c}^{1}(\Omega)$, we have
\begin{equation}\label{wksol}
\int_{\Omega}\mathcal{B}(x,\nabla u)\cdot\nabla\phi(x)\,dx+\int_{\mathbb{R}^N}\int_{\mathbb{R}^N}\mathcal{A}u(x,y)(\phi(x)-\phi(y))\,d\mu
=\int_{\Omega}f(x){u^{-\gamma(x)}}\phi(x)\,dx,
\end{equation}
where $\mathcal{A}u(x,y)$ and $d\mu$ are as in \eqref{not}.
\end{Definition}

\begin{Remark}\label{defrmk2}
Lemma \ref{defineq} and Lemma \ref{locnon1} ensure that  Definition \ref{psintest} is well stated.
\end{Remark}

\section{Regularity results}
The main challenge in singular problems is that the classical regularity results cannot be directly applied. One approach to overcome such a difficulty is the method of approximation as introduced in Boccardo-Orsina \cite{BocOr} for a singular Laplace equation. We employ a similar approximation technique for a mixed anisotropic problem. 
We make use of Carmona-Mart\'{\i}nez-Aparicio \cite{CMP} to accommodate the presence of variable singularity. 
Moreover, we adopt the technique from Garain-Ukhlov \cite{GU} for mixed singular problems. One of the main ingredients used in the approximation technique is the regularity theory of the associated non-singular problem.

Regularity results are well studied in both the local  and nonlocal cases, see \cite{BLS, Di} and the references therein. 
In the mixed linear case several regularity estimates, for example Harnack's inequality, are proved by probabilistic methods in Foondun \cite{Fo} and Chen-Kim-Song-Vondra\v{c}ek \cite{CKSV}. 
Recently, Biagi-Dipierro-Salort-Valdinoci-Vecchi \cite{Biagi2, BDVV, SV} obtained various regularity results, including many qualitative properties of solutions, by a purely analytic approach.
Such results are very useful in the mixed singular problems in the linear case in \cite{AroRad}. In the quasilinear case, for the mixed local and nonlocal $p$-Laplace equation, 
an analytic approach to regularity theory is discussed in \cite{GK}. This can be used to develop the existence theory for mixed singular problems in the quasilinear setting. 
Recently, there has been a lot of interest in mixed problems, see \cite{BDVV2, Biagi1, BVDV, DeMin, GL} for the elliptic case and \cite{FSZ, GK2, GKwh, SZ} for the parabolic case. 
Regularity results for the mixed parabolic problems including the anisotropic case can be found in \cite{GK2}.

In this section, we discuss the necessary regularity results for \eqref{genmaineqn}.
As far as we are aware, the regularity theory for elliptic mixed problems in the presence of the anisotropic operator has not been studied to date.
However, following De-Giorgi theory adopted in \cite{GK}, it is possible to obtain local boundedness, Harnack inequality, weak Harnack inequality for weak solutions of \eqref{genmaineqn}. 
The energy estimates below will be important for us. We begin with an energy estimate for weak subsolutions of \eqref{genmaineqn}.

\begin{Lemma}\label{energyest}
Let $u$ be a weak subsolution of \eqref{genmaineqn} and denote $w=(u-k)^{+}$ with $k\in\mathbb{R}$. 
There exists a constant $c=c(p,\Lambda,C_1,C_2)$ such that 
\begin{equation}\label{energyesteqn}
\begin{split}
&\int_{B_r(x_0)}\psi(x)^p|\nabla w(x)|^p\,dx+\int_{B_r(x_0)}\int_{B_r(x_0)}|w(x)\psi(x)-w(y)\psi(y)|^p \,d\mu\\
&\leq C\bigg(\int_{B_r(x_0)}w(x)^p |\nabla\psi(x)|^p\,dx+\int_{B_r(x_0)}\int_{B_r(x_0)}{\max\{w(x),w(y)\}^p|\psi(x)-\psi(y)|^p}\,d\mu\\
&\qquad+\esssup_{x\in\supp\psi}\int_{{\mathbb{R}^n\setminus B_r(x_0)}}{\frac{w(y)^{p-1}}{|x-y|^{n+ps}}}\,dy
\cdot\int_{B_r(x_0)}w\psi^p\,dx\bigg),
\end{split}
\end{equation}
whenever $B_r(x_0)\subset\Omega$ and $\psi\in C_c^{\infty}(B_r{(x_0)})$ is a nonnegative function.
\end{Lemma}

\begin{proof}
Let $u$ be a weak subsolution of \eqref{genmaineqn}. For $w=(u-k)^+$, by choosing $\phi=w\psi^p$ as a test function in \eqref{weaksubsupsoln}, we obtain
\begin{equation}\label{genenergytest}
\begin{split}
0&\geq\int_{B_r(x_0)}\mathcal{B}(x,u,\nabla u)\cdot\nabla(w(x)\psi(x)^p)\,dx\\
&\qquad+\int_{\mathbb{R}^n}\int_{\mathbb{R}^n}\mathcal{A}(u(x,y))(w(x)\psi(x)^p-w(y)\psi(y)^p)\,d\mu\\
&=I+J.
\end{split}
\end{equation}
As in \cite[Page 14, Proposition 3.1]{Verenacontinuity}, using Young's inequality and {\rm(H1)}, we obtain
\begin{equation}\label{genestI}
\begin{split}
I&\geq c\int_{B_r(x_0)}\psi(x)^p|\nabla w(x)|^p\,dx-C\int_{B_r(x_0)}w(x)^p|\nabla\psi(x)|^p\,dx,
\end{split}
\end{equation}
for some positive constants $c,C$ depending on $p,C_1,C_2$. The rest of the proof follows along the lines of the proof of \cite[Lemma 3.1]{GK}.
\end{proof}

\begin{Remark}\label{genrmkeng}
If $u$ is a weak supersolution of \eqref{genmaineqn}, the estimate in \eqref{energyesteqn} holds by replacing $w$ with $(u-k)^{-}$. This can be proved in a similar way as Lemma \ref{energyest}.
\end{Remark}

A tail quantity will appear in regularity estimates for nonlocal problems.

\begin{Definition}\label{def.tail}
Let $u$ be a weak subsolution, or a weak supersolution, of \eqref{genmaineqn} as in Definition \ref{gensubsupsolution}.
The tail of $u$ with respect to a ball $B_r(x_0)$ is defined by
\begin{equation}\label{loctail}
\Tail(u;x_0,r)=\bigg(r^{p}\int_{\mathbb{R}^N\setminus B_r(x_0)}\frac{|u(y)|^{p-1}}{|y-x_0|^{N+ps}}\,dy\bigg)^\frac{1}{p-1}.
\end{equation}
\end{Definition}

The following tail estimate will be useful for us.

\begin{Lemma}\label{Tail}
Let $u$ be a weak solution of \eqref{genmaineqn} such that $u\geq 0$ in $B_R(x_0)\subset\Omega$.
Then there exists a constant $c=c(n,p,s,\Lambda,C_1,C_2)$ such that
\begin{equation}\label{tailest}
\Tail(u^{+};x_0,r)\leq c\esssup_{x\in B_r(x_0)}u(x)+c\Big(\frac{r}{R}\Big)^\frac{p}{p-1}\Tail(u^{-};x_0,R),  
\end{equation}
whenever $0<r<R$ with $r\in(0,1]$.
Here $\Tail(\cdot)$ is given by \eqref{loctail}.
\end{Lemma}

\begin{proof}
Let $M=\esssup_{B_r(x_0)}\,u$ and $\psi\in C_c^{\infty}(B_r(x_0))$ be a cutoff function such that
$0\leq\psi\leq 1$ in $B_r(x_0)$, $\psi=1$ in $B_{\frac{r}{2}}(x_0)$ and $|\nabla\psi|\leq\frac{8}{r}$ in $B_r(x_0)$.
By letting $w=u-2M$. We observe that
\begin{equation}
\begin{split}
\mathcal{B}(x,u,\nabla u)\cdot\nabla(w(x)\psi(x)^p)
&=\mathcal{B}(x,u,\nabla w)\cdot\nabla(w(x)\psi(x)^p)\\
&\geq C_1\psi(x)^p|\nabla w(x)|^p+p\psi(x)^{p-1}w(x)\mathcal{B}(x,u,\nabla w)\nabla\psi(x)\\
&\geq \frac{C_1}{2}\psi(x)^p|\nabla w(x)|^p-C(C_1,C_2,p)w(x)^p|\nabla\psi(x)|^p,
\end{split}
\end{equation}
where $C_1,C_2$ are given by {\rm(H1)}.
The result follows by choosing $\phi=w\psi^p$ as a test function in \eqref{weaksubsupsoln} and proceeding as in the proof of \cite[Lemma 6.1]{GK}.
\end{proof}

Next, we prove the following energy estimate for weak supersolutions of \eqref{genmaineqn}.

\begin{Lemma}\label{energyforrev}
Let $1<q<p$ and $d>0$.
Assume that $u$ is a weak supersolution of \eqref{genmaineqn} such that $u\geq 0$ in $B_R(x_0)\subset\Omega$ and let $w=(u+d)^\frac{p-q}{p}$. 
Then there exists a constant $c=c(p,\Lambda,C_1,C_2)$ such that
\begin{equation}\label{energyrevest1}
\begin{split}
&\int_{B_r(x_0)}\psi(x)^p|\nabla w(x)|^p\,dx
\leq c\bigg(\frac{(p-q)^p}{(q-1)^\frac{p}{p-1}}\int_{B_r(x_0)}w^p|\nabla\psi|^p\,dx\\
&\qquad+\frac{(p-q)^p}{(q-1)^{p}}\int_{B_r(x_0)}\int_{B_r(x_0)}\max\{w(x),w(y)\}^p|\psi(x)-\psi(y)|^p\,d\mu\\
&\qquad+\frac{(p-q)^p}{(q-1)}\bigg(\esssup_{z\in\supp\psi}\int_{\mathbb{R}^n\setminus B_r(x_0)}K(z,y)\,dy
+d^{1-p}R^{-p}\Tail(u_{-};x_0,R)^{p-1}\bigg)\\
&\qquad\qquad\cdot\int_{B_r(x_0)}w(x)^p\psi(x)^p\,dx\bigg),
\end{split}
\end{equation}
whenever $B_r(x_0)\subset B_{\frac{3R}{4}}(x_0)$ and $\psi\in C_c^{\infty}(B_r(x_0))$ is a nonnegative function.
Here $\Tail(\cdot)$ is defined in \eqref{tailest}.
\end{Lemma}

\begin{proof}
Let $d>0$, $v=u+d$ and $1+\epsilon\le q\le p-\epsilon$ for $\epsilon>0$ small enough. By choosing $\phi=v^{1-q}\psi^p$ as a test function in \eqref{weaksubsupsoln}, we obtain
\begin{equation}\label{revestfinal}
\begin{split}
0&\leq \int_{B_r(x_0)}\mathcal{B}(x,u,\nabla u)\cdot\nabla( v(x)^{1-q}\psi(x)^p)\,dx\\
&\qquad+\int_{B_r(x_0)}\int_{B_r(x_0)}\mathcal{A}v(x,y)(v(x)^{1-q}\psi(x)^p- v(y)^{1-q}\psi(y)^p)\,d\mu\\
&\qquad+2\int_{\mathbb{R}^n\setminus B_r(x_0)}\int_{B_r(x_0)}\mathcal{A}( v(x,y)) v(x)^{1-q}\psi(x)^p\,d\mu.
\end{split}
\end{equation}
By applying {\rm(H1)}, we have
\begin{equation}\label{genI_1est}
\begin{split}
I_1&=\int_{B_r(x_0)}\mathcal{B}(x,u,\nabla u)\cdot\nabla(v(x)^{1-q}\psi(x)^p)\,dx\\
&=\int_{B_r(x_0)}\mathcal{B}(x,u,\nabla u)\cdot\big((1-q)v(x)^{-q}\psi(x)^p\nabla v(x)+p\psi(x)^{p-1}v(x)^{1-q}\nabla\psi(x)\big)\,dx\\
&\leq C_2(1-q)\int_{B_r(x_0)}v(x)^{-q}|\nabla v(x)|^p\psi(x)^p\,dx\\
&\qquad+pC_2\int_{B_r(x_0)}v^{1-q}|\nabla\psi(x)||\nabla v(x)|^{p-1}\psi(x)^{p-1}\,dx.
\end{split}
\end{equation}
The rest of the proof follows as in the proof of \cite[Lemma 3.3]{GK}.
\end{proof}

The next result shows that expansion of positivity holds for the general mixed problem \eqref{genmaineqn}.
\begin{Lemma}\label{DGLemma}
Let $u$ be a weak supersolution of \eqref{genmaineqn} such that $u\geq 0$ in $B_R(x_0)\subset\Omega$. 
Assume $k\geq 0$ and there exists $\tau\in(0,1]$ such that
\begin{equation}\label{expangiven}
\big|B_r(x_0)\cap\{u\geq k\}\big|\geq \tau|B_r(x_0)|,
\end{equation}
for some $r\in(0,1]$ with $0<r<\frac{R}{16}$. There exists a constant $\delta=\delta(n,p,s,\Lambda,\tau)\in(0,\frac{1}{4})$ such that
\begin{equation}\label{expan}
\essinf_{x\in B_{4r}(x_0)}u(x)\geq\delta k-\Big(\frac{r}{R}\Big)^\frac{p}{p-1}\Tail(u^{-};x_0,R),
\end{equation} 
where $\Tail(\cdot)$ is given by \eqref{loctail}.
\end{Lemma}

\begin{proof}
The proof is similar to the proof of \cite[Lemma 7.1]{GK}.
\end{proof}

The following preliminary version of a weak Harnack inequality follows from Lemma \ref{DGLemma} and the proof of \cite[Lemma 8.1]{GK}, 

\begin{Lemma}\label{WeakHarnacklemma}
Let $u$ be a weak supersolution of \eqref{genmaineqn} such that $u\geq 0$ in $B_R(x_0)\subset\Omega$. 
Then there exist constants $\eta=\eta(n,p,s,\Lambda)\in(0,1)$ and $c=c(n,p,s,\Lambda)\geq 1$ such that
\begin{equation}\label{wk}
\bigg(\fint_{B_r(x_0)}u(x)^{\eta}\,dx\bigg)^\frac{1}{\eta}
\leq c\essinf_{x\in B_r(x_0)}u(x)+c\Big(\frac{r}{R}\Big)^\frac{p}{p-1}\Tail(u^{-};x_0,R),
\end{equation}
whenever $B_r(x_0)\subset B_R(x_0)$ with $r\in(0,1]$.
Here $\Tail(\cdot)$ is defined in \eqref{loctail}.
\end{Lemma}

The next result asserts that  weak subsolutions are locally bounded.
This can be proved as in the proof of \cite[Theorem 4.2]{GK} by taking into account Lemma \ref{energyest}.

\begin{Theorem}\label{locbdd}(\textbf{Local boundedness}).
Suppose $\mathcal{B}(x,u,\nabla u)$ satisfies the hypothesis {\rm(H1)}. Let $u$ be a weak subsolution of \eqref{genmaineqn}. 
Then there exists a positive constant $c=c(n,p,s,\Lambda,C_1,C_2)$, such that
\begin{equation}\label{locbd}
\esssup_{x\in B_{\frac{r}{2}}(x_0)}u(x)
\leq \delta \Tail\big(u^{+};x_0,\tfrac{r}{2}\big)+c\delta^{-\frac{(p-1)\kappa}{p(\kappa-1)}}\bigg(\fint_{B_r(x_0)}u^{+}(x)^p\,dx\bigg)^\frac{1}{p},
\end{equation}
whenever $B_r(x_0)\subset\Omega$ with $r\in(0,1]$ and $\delta\in(0,1]$. Here $\Tail(\cdot)$ is as defined in \eqref{loctail} and
\begin{equation}\label{kappa}
\kappa=
\begin{cases}
\frac{N}{N-p},&\text{if}\quad 1<p<N,\\
2,&\text{if}\quad p\geq N.
\end{cases}
\end{equation} 
\end{Theorem}

A Harnack inequality follows as the proof of \cite[Theorem 8.3]{GK} but
taking into account Lemma \ref{locbdd}, Lemma \ref{tailest} and Lemma \ref{WeakHarnacklemma}.

\begin{Theorem}\label{regthm2}(\textbf{Harnack inequality}).
Suppose $\mathcal{B}(x,u,\nabla u)$ satisfies {\rm(H1)}. Let $u$ be a weak solution of \eqref{genmaineqn} such that $u\geq 0$ in $B_R(x_0)\subset\Omega$. 
Then there exists a constant $c=c(n,p,s,\Lambda,C_1,C_2)$ such that
\begin{equation}\label{Harnackest}
\esssup_{x\in B_{\frac{r}{2}}(x_0)}u(x)
\leq c\essinf_{x\in B_r(x_0)}u(x)+c\Big(\frac{r}{R}\Big)^\frac{p}{p-1}\Tail(u^{-};x_0,R),
\end{equation}
whenever $B_r(x_0)\subset B_\frac{R}{2}(x_0)$ and $r\in(0,1]$.
Here $\Tail(\cdot)$ is given by \eqref{loctail}.
\end{Theorem}

The proof of the next result follows by taking into account Lemma \ref{energyforrev}, Lemma \ref{WeakHarnacklemma} and proceeding along the lines of the proof of \cite[Theorem 8.4]{GK}.

\begin{Theorem}\label{regthm3}(\textbf{Weak Harnack inequality}).
Suppose $\mathcal{B}(x,u,\nabla u)$ satisfies {\rm(H1)}. Let $u$ be a weak supersolution of \eqref{genmaineqn} such that $u\geq 0$ in $B_R(x_0)\subset\Omega$. 
Then there exists a constant $c=c(n,p,s,\Lambda,C_1,C_2)$ such that
\begin{equation}\label{wkest}
\bigg(\fint_{B_{\frac{r}{2}}(x_0)}u(x)^{l}\,dx\bigg)^\frac{1}{l}
\leq c\essinf_{x\in B_{r}(x_0)}u(x)+c\Big(\frac{r}{R}\Big)^\frac{p}{p-1}\Tail(u^{-};x_0,R),
\end{equation}
whenever $B_r(x_0)\subset B_{\frac{R}{2}}(x_0)$, $r\in(0,1]$ and  $0<l<\kappa (p-1)$.
Here $\kappa$ and $\Tail(\cdot)$ are defined in \eqref{kappa} and \eqref{loctail} respectively.
\end{Theorem}

\section{Existence results}

\subsection{Statement of the existence results}
The main existence results of the paper are stated below. To this end, we recall the definition of
$$
\Omega_{\delta}=\{x\in\Om:\text{dist}(x,\partial\Om)<\delta\},\quad\delta>0.
$$

\begin{Theorem}\label{thm2}
Let $1<p<\infty$ and $\gamma\in C(\overline{\Omega})$ with $\gamma(x)>0$ for every $x\in\overline{\Omega}$. Suppose $\mathcal{B}(x,\nabla u)$ satisfies {\rm(H1)} and {\rm(H2)}.
\begin{enumerate}
\item[{\rm(a)}] Assume that there exists $\delta>0$ such that $0<\gamma(x)\leq 1$ for every in $x\in\Om_{\delta}$ and $f\in L^m(\Om)\setminus\{0\}$ is a nonnegative function with
\begin{equation}\label{m1}
m=
\begin{cases}
(p^*)',\text{ if }1<p<N,\\
>1,\text{ if }p=N,\\
1,\text{ if }p>N.
\end{cases}
\end{equation}
Then \eqref{sinmeqn1} admits a unique weak solution $u\in W_0^{1,p}(\Om)$.
\item[{\rm(b)}] Assume that there exists $\delta>0$ and $\gamma^*>1$ such that $\|\gamma\|_{L^\infty(\Om_\delta)}\leq\gamma^*$ and that $f\in L^m(\Om)\setminus\{0\}$ is a nonnegative function with
\begin{equation}\label{m3}
m=
\begin{cases}
\Big(\frac{(\gamma^{*}+p-1){p^*}}{p\gamma^{*}}\Big)',\text{ if } 1< p<N,\\ 
\Big(\frac{(\gamma^{*}+p-1)l}{p\gamma^{*}}\Big)',\text{ if } p=N \text{ and } \frac{p\gamma^{*}}{\gamma^{*}+p-1}<l<\infty,\\
1, \text{ if } p>N.
\end{cases}
\end{equation}
Then  \eqref{sinmeqn1} admits a weak solution $u\in W^{1,p}_{\mathrm{loc}}(\Om)$ such that $u^\frac{\gamma^*+p-1}{p}\in W_0^{1,p}(\Om)$.
\end{enumerate}
\end{Theorem}

\begin{Remark}\label{rmknew}
Under the hypothesis that $\gamma\in C^1(\overline{\Omega})$ and $1<p<N$, Theorem \ref{thm2}-(a) generalizes \cite[Theorem 4.3]{Biroud} and Theorem \ref{thm2}-(b) generalizes \cite[Theorem 4.4]{Biroud}.
\end{Remark}

Furthermore, when $\gamma(x)=\gamma$ is a positive constant, we can extend the class of functions $f$ in the existence results as stated in Theorem \ref{cthm1} and Theorem \ref{cthm3}.
\begin{Theorem}\label{cthm1}
Let $1<p<\infty$ and $0<\gamma<1$. Suppose $\mathcal{B}(x,\nabla u)$ satisfies {\rm(H1)} and {\rm(H2)}. Assume that $f\in L^m(\Omega)\setminus\{0\}$ is a nonnegative function with
\begin{equation}\label{cm}
    m=
\begin{cases}
    \big(\frac{p^{*}}{1-\gamma}\big)',& \text{if } 1< p<N, \\
   >1,& \text{if } p=N, \\
    1,& \text{if } p>N.
\end{cases}
\end{equation}
Then \eqref{sinmeqn1} admits a weak solution $u\in W_0^{1,p}(\Omega)$.
\end{Theorem}

\begin{Theorem}\label{cthm2}
Let $1<p<\infty$ and $\gamma=1$. Suppose $\mathcal{B}(x,\nabla u)$ satisfies {\rm(H1)} and {\rm(H2)}. 
Then \eqref{sinmeqn1} admits a weak solution $u\in W_0^{1,p}(\Om)$ for any nonnegative $f\in L^1(\Om)\setminus\{0\}$.
\end{Theorem}

\begin{Theorem}\label{cthm3}
Let $1<p<\infty$ and $\gamma>1$. Suppose $\mathcal{B}(x,\nabla u)$ satisfies {\rm(H1)} and {\rm(H2)}. 
Then  \eqref{sinmeqn1} admits a weak solution $u\in W^{1,p}_{\mathrm{loc}}(\Om)$ such that $u^\frac{\delta+p-1}{p}\in W_0^{1,p}(\Om)$ for any nonnegative $f\in L^1(\Om)\setminus\{0\}$.
\end{Theorem}

In order to prove the existence results above, we follow the approximation method in \cite{BocOr}. To this end, we investigate the approximating problem
\begin{equation}\label{approx}
-\mathcal{B}(x,\nabla u)+\mathcal{L}(u)=f_n(x)\Big(u^{+}(x)+\frac{1}{n}\Big)^{-\gamma(x)}\text{ in }\Omega,\,u=0\text{ on }\partial\Omega,
\end{equation}
where $f_n(x)=\min\{f(x),n\}$ for $x\in\Omega$, $n\in\mathbb{N}$, $f\in L^1(\Omega)\setminus\{0\}$ is a nonnegative function, $\gamma\in C(\overline{\Omega})$ is positive, $\mathcal{B}(x,\nabla u)$ satisfy {\rm(H1)} and {\rm(H2)} and $1<p<\infty$.


\begin{Lemma}\label{exisapprox}
\begin{enumerate}
\item[{\rm(a)}]
For every $n\in\mathbb{N}$, the problem \eqref{approx} admits a positive weak solution $u_{n}\in W_0^{1,p}(\Om)\cap L^\infty(\Om)$.
\item[{\rm(b)}]
$u_{n+1}\geq u_n$ for every $n\in\mathbb{N}$ and $u_n$ is unique for every $n\in\mathbb{N}$.
\item[{\rm(c)}]
For every $\Omega'\Subset\Omega$, there exists a positive constant $C(\Omega')$ such that $u_n\geq C(\Omega')>0$ in $\Omega'$.
\end{enumerate}
\end{Lemma}

\begin{proof}
\textbf{Existence:} Let $V=W_0^{1,p}(\Om)$ with the norm $\|u\|=\|\nabla u\|_{L^p(\Om)}$ and let $V'$ be the dual of $V$. 
For every $n\in\mathbb{N}$ we define $T:V\to V'$ by
\begin{align*}
\langle T(v),\phi\rangle&=\int_{\Om}\mathcal{B}(x,\nabla v)\cdot\nabla\phi(x)\,dx
+\int_{\mathbb{R}^N}\int_{\mathbb{R}^N}\mathcal{A}(v(x,y))(\phi(x)-\phi(y))\,d\mu\\
&\qquad-\int_{\Om}f_n(x)\Big(v^{+}(x)+\frac{1}{n}\Big)^{-\gamma(x)}\phi(x)\,dx
\end{align*}
for every $v,\phi\in V$.
We observe that $T$ is well defined. Using {\rm(H1)}, H$\ddot{\text{o}}$lder's inequality, Lemma \ref{locnon1} and Lemma \ref{emb}, we obtain 
\begin{align*}
|\langle T(v),\phi\rangle\big|
&\leq C_2\int_{\Om}|\nabla v(x)|^{p-1}|\nabla\phi(x)|\,dx
+\Lambda[v]_{s,p}^{p-1}[\phi]_{s,p}+n^{1+\|\gamma\|_{L^\infty(\Om)}}\int_{\Om}|\phi(x)|\,dx\\
&\leq C(\|v\|^{p-1}\|\phi\|+\|\phi\|)\leq C\|\phi\|,
\end{align*}
for some constant $C$. 

\textbf{Coercivity:} By H\"older's inequality,  {\rm(H1)} and Lemma \ref{emb}, we obtain
\begin{align*}
\langle T(v),v\rangle&=\int_{\Om}\mathcal{B}(x,\nabla v)\cdot\nabla v(x)\,dx+\int_{\mathbb{R}^N}\int_{\mathbb{R}^N}\mathcal{A}(v(x,y))(v(x)-v(y))\,d\mu\\
&\qquad-\int_{\Om}f_n(x)\Big(v^+(x) +\frac{1}{n}\Big)^{-\gamma(x)}v(x)\,dx\\
&\geq C_1\|v\|^{p}-C\|v\|.
\end{align*}
Since $1<p<\infty$, we may conclude that $T$ is coercive.

\textbf{Demicontinuity:} Let $v_k\to v$ in the norm of $V$ as $k\to\infty$. Then $\nabla v_k\to \nabla v$ in $L^p(\Omega)$ as $k\to\infty$. Therefore up to a subsequence, still denoted by $v_k$, we have $v_k(x)\to v(x)$ and $\nabla v_{k}(x)\to\nabla v(x)$ pointwise for almost every $x\in\Omega$. Since the function $\mathcal{B}(x,\cdot)$ is continuous in the second variable, we have 
$\mathcal{B}(x,\nabla v_{k})\to \mathcal{B}(x,\nabla v)$
pointwise for almost every $x\in\Omega$ as $k\to\infty$. 
By {\rm(H1)} and the norm boundedness of the sequence $(v_k)_{k\in\mathbb N}$, we obtain
\[
||\mathcal{B}(x,\nabla v_{k})||^{p'}_{L^{p'}(\Omega)}
=\int_{\Omega}|\mathcal{B}(x,\nabla v_{k})|^{p'}\,dx
\leq C_2^{p'}\int_{\Omega}|\nabla v_{k}(x)|^p\,dx\leq C,
\]
for some positive constant $C$, independent of $k$. Therefore, up to a subsequence
$$
\mathcal{B}(x,\nabla v_{k}(x))\rightharpoonup \mathcal{B}(x,\nabla v(x))
$$ 
weakly in $L^{p'}(\Om)$ as $k\to\infty$ and since the weak limit is independent of the choice of the subsequence, the above weak convergence holds for the original sequence. 
Since $\phi\in V$ implies that $\nabla\phi\in L^p(\Omega)$, by the weak convergence, we obtain
\begin{equation}\label{d1}
\lim_{k\to\infty}\int_{\Om}\mathcal{B}(x,\nabla v_k)\cdot\nabla\phi(x)\,dx
=\int_{\Om}\mathcal{B}(x,\nabla v)\cdot\nabla\phi(x)\,dx.
\end{equation}
Since $(v_k)_{k\in\mathbb N}$ is a bounded bounded sequence in $V$, Lemma \ref{locnon1} implies that
$(\mathcal{A}v_kK^\frac{1}{p'})_{k\in\mathbb N}$
is a bounded sequence in $L^{p'}(\mathbb{R}^N\times\mathbb{R}^N)$.
Since $\phi\in V$, we have
\[
(\phi(x)-\phi(y))K(x,y)^\frac{1}{p}\in L^p(\mathbb{R}^N\times\mathbb{R}^N).
\]
Therefore, by the weak convergence, we have
\begin{equation}\label{nonloclimnew}
\lim_{k\to\infty}\int_{\mathbb{R}^N}\int_{\mathbb{R}^N}\mathcal{A}v_k(x,y)(\phi(x)-\phi(y))\,d\mu
=\int_{\mathbb{R}^N}\int_{\mathbb{R}^N}\mathcal{A}v(x,y)(\phi(x)-\phi(y))\,d\mu.
\end{equation}
Moreover, we have $v_k(x)\to v(x)$ pointwise almost everywhere in $\Om$ and
$$
\left|f_n\Big(v_k^{+}(x)+\frac{1}{n}\Big)^{-\gamma(x)}\phi(x)-f_n(x)\Big(v^{+}(x)+\frac{1}{n}\Big)^{-\gamma(x)}\phi(x)\right|
\leq n^{1+\|\gamma\|_{L^\infty(\Om)}}|\phi(x)|
$$
for every $\phi\in V$.
By the dominated convergence theorem, we have
\begin{equation}\label{d2}
\lim_{k\to\infty}\int_{\Om}f_n(x)\Big(v_k^{+}(x)+\frac{1}{n}\Big)^{-\gamma(x)}\phi(x)\,dx
=\int_{\Om}f_n(x)\Big(v^{+}(x)+\frac{1}{n}\Big)^{-\gamma(x)}\phi(x)\,dx
\end{equation}
for every $\phi\in V$ and $n\in\mathbb{N}$.
Therefore, from \eqref{d1}, \eqref{nonloclimnew} and \eqref{d2}, it follows that
$$
\lim_{k\to\infty}\langle T(v_k),\phi\rangle=\langle T(v),\phi\rangle
$$
for every $\phi\in V$
and hence $T$ is demicontinuous.

\textbf{Monotonicity:} Let $v_1,v_2\in V$. We observe that
\begin{equation}\label{tmn}
\begin{split}
&\langle T(v_1)-T(v_2),v_1 -v_2\rangle
=\int_{\Om}(\mathcal{B}(x,\nabla v_1)-\mathcal{B}(x,\nabla v_2))\cdot(\nabla v_1(x)-\nabla v_2(x))\,dx\\
&\qquad+\int_{\mathbb{R}^N}\int_{\mathbb{R}^N}\big(\mathcal{A}v_1(x,y)-\mathcal{A}v_2(x,y)\big)\big((v_1(x)-v_2(x))-(v_1(y)-v_2(y))\big)\,d\mu\\
&\qquad-\int_{\Om}f_n(x)\left(\Big(v_1^{+}(x)+\frac{1}{n}\Big)^{-\gamma(x)}-\Big(v_2^{+}(x)+\frac{1}{n}\Big)^{-\gamma(x)}\right)(v_1(x)- v_2(x))\,dx.
\end{split}
\end{equation}
Note that by {\rm(H2)} and Lemma \ref{alg}, the first and second integrals above are nonnegative. Moreover, it is easy to see that the third integral above is nonpositive. 
This shows that $T$ is monotone.

\textbf{Existence:}
Theorem \ref{MB} implies that $T$ is surjective. Hence for every $n\in\mathbb{N}$, there exists $u_n\in V$ such that
\begin{equation}\label{vneqn}
\begin{split}
\int_{\Om}\mathcal{B}(x,\nabla u_n)\cdot\nabla\phi(x)&\,dx
+\int_{\mathbb{R}^N}\int_{\mathbb{R}^N}\mathcal{A}u_n(x,y)(\phi(x)-\phi(y))\,d\mu\\
&=\int_{\Om}f_n(x)\Big(u_n^{+}(x) +\frac{1}{n}\Big)^{-\gamma(x)}\,\phi(x)\,dx
\end{split}
\end{equation}
for every $\phi\in V$.
This shows that \eqref{approx} has a  weak solution $u_{n}\in W_0^{1,p}(\Om)$ for every $n\in\mathbb{N}$.

\textbf{Boundedness:} Proceeding along the lines of the proof of \cite[Lemma 3.1]{GU}, it follows that $u_n\in L^\infty(\Om)$ for every $n\in\mathbb{N}$.

\textbf{Monotonicity and uniqueness of $u_n$:} Choosing $\phi = (u_{n}-u_{n+1})^+$ as a test function in \eqref{vneqn}, we have
\begin{equation}\label{monopf}
\begin{split}
&\int_{\Omega}(\mathcal{B}(x,\nabla u_n)-\mathcal{B}(x,\nabla u_{n+1}))\cdot\nabla(u_n(x)-u_{n+1}(x))^{+}\,dx\\
&\qquad\qquad+\int_{\mathbb{R}^N}\int_{\mathbb{R}^N}   (  \mathcal{A}u_n(x,y)   -\mathcal{A}(u_{n+1}(x,y)) )
   \big((u_n(x)-u_{n+1(x)}\big)^+\\
   &\qquad\qquad\qquad\qquad\qquad-\big(u_n(y)-u_{n+1}(y))^+\big)\,d\mu\\
&\qquad=\int_{\Omega}\left(f_{n}(x)\Big(u_{n}(x)+\frac{1}{n}\Big)^{-\gamma(x)}-f_{n+1}(x)\Big(u_{n+1}(x)+\frac{1}{n+1}\Big)^{-\gamma(x)}\right)\\
&\qquad\qquad\qquad\qquad\qquad\cdot(u_n(x)-u_{n+1}(x))^{+}\,dx\\
&\qquad\leq\int_{\Omega}f_{n+1}(x)\left(\Big(u_{n}(x)+\frac{1}{n}\Big)^{-\gamma(x)}-\Big(u_{n+1}(x)+\frac{1}{n+1}\Big)^{-\gamma(x)}\right)\\
&\qquad\qquad\qquad\qquad\qquad\cdot(u_n(x)-u_{n+1}(x))^{+}\,dx\leq 0,
\end{split}
\end{equation}
where the last inequality above follows using $f_{n}(x) \leq f_{n+1}(x)$ and the positivity of $\gamma$. Following the same arguments from the proof of \cite[Lemma $9$]{Ling}, we obtain
\begin{equation}\label{monnon}
\int_{\mathbb{R}^N}\int_{\mathbb{R}^N}\big(\mathcal{A}u_n(x,y)-\mathcal{A}u_{n+1}(x,y)\big)\big((u_n(x)-u_{n+1}(x))^+-(u_n(y)-u_{n+1}(y))^+\big)\,d\mu\geq 0.
\end{equation}
By \eqref{monnon} in \eqref{monopf}, we arrive at
$$
\int_{\Omega}(\mathcal{B}(x,\nabla u_n)-\mathcal{B}(x,\nabla u_{n+1}))\cdot\nabla(u_n(x)-u_{n+1}(x))^{+}\,dx\leq 0.
$$
Thus, by {\rm(H2)}, we may conclude that $\nabla u_n=\nabla u_{n+1}$ in $\{x\in\Om:u_n(x)>u_{n+1}(x)\}$. Hence, we have $u_{n+1}\geq u_n$ in $\Om$. 
The uniqueness of $u_n$ follows similarly. 

\textbf{Positivity:} Let 
$$
g_n(x)=f_n(x)\Big(u_n^{+}(x) +\frac{1}{n}\Big)^{-\gamma(x)}.
$$
Choosing $\phi=\min\{u_n,0\}=(u_n)_-$ as a test function in \eqref{vneqn} and using the nonnegativity of the right-hand side of \eqref{vneqn} and {\rm(H1)}, we have
\begin{equation}\label{posi}
\begin{split}
C_1\int_{\Omega}|\nabla (u_n(x)_-)|^p&\,dx+\int_{\mathbb{R}^N}\int_{\mathbb{R}^N}\mathcal{A}u_n(x,y)(u_n(x)_--u_n(y)_-)\,d\mu\\
&=\int_{\Omega}g_n(x)u_n(x)_-\,dx\leq 0,
\end{split}
\end{equation}
Proceeding similarly as in the proof of the estimate \cite[(3.13)]{GU}, the second integral above is nonnegative. Hence, we get
$$
\int_{\Omega}|\nabla (u_n(x)_-)|^p\,dx=0,
$$
which gives $u_n\geq 0$ in $\Om$. Since $g_n\neq 0$ in $\Om$, by Theorem \ref{regthm3} for every
$\Omega'\Subset\Om$, there exists a positive constant $C=C(\Omega',n)$ such that 
$u_n\geq C(\Omega',n)>0$ in $\Omega'$.
Using the monotonicity property, we get $u_n\geq u_1$ in $\Om$ and thus, for every $\Omega'\Subset\Om$, there exists a positive constant $C(\Omega')$ (independent of $n$) such that $u_n\geq C(\Omega')>0$ in $\Omega'$.
\end{proof}

Next we consider a priori estimates for the solution of the approximating problem. 

\begin{Lemma}\label{apun}
Let $(u_n)_{n\in\mathbb{N}}$ be the sequence of solutions of \eqref{approx} given by Lemma \ref{exisapprox}.
\begin{enumerate}
\item[{\rm(a)}] Assume that there exists $\delta>0$ such that $0<\gamma(x)\leq 1$ for every $x\in\Omega_\delta$ and $f\in L^m(\Omega)\setminus\{0\}$ is a nonnegative function, where $m$ is as in \eqref{m1}.
Then $(u_n)_{n\in\mathbb{N}}$ is a bounded sequence in $W_0^{1,p}(\Omega)$.
\item[{\rm(b)}] Assume that $\gamma(x)=\gamma$, where $0<\gamma<1$ is a constant and $f\in L^m(\Omega)\setminus\{0\}$ is a nonnegative function, where $m$ is defined in \eqref{cm}.
Then $(u_n)_{n\in\mathbb{N}}$ is a bounded sequence in $W_0^{1,p}(\Om)$
\item[{\rm(c)}]  Assume that $\gamma(x)=1$ and  $f\in L^1(\Omega)\setminus\{0\}$ is a nonnegative function.
Then $(u_n)_{n\in\mathbb{N}}$ is a bounded sequence in $W_0^{1,p}(\Om)$.
\item[{\rm(d)}] Assume that there exists $\delta>0$ and $\gamma^*>1$ such that $\|\gamma\|_{L^\infty(\Omega_\delta)}\leq\gamma^*$ and $f\in L^m(\Omega)\setminus\{0\}$ is a nonnegative function, where $m$ is as in \eqref{m3}.
Then $(u_n)_{n\in\mathbb{N}}$ and  $(u_n^{(\gamma^{*}+p-1)/p})_{n\in\mathbb{N}}$ are bounded sequences in $W_{\mathrm{loc}}^{1,p}(\Omega)$ and $W_0^{1,p}(\Omega)$, respectively.
\item[{\rm(e)}] Assume that $\gamma(x)=\gamma$, where $\gamma>1$ is a constant and $f\in L^1(\Omega)\setminus\{0\}$ is a nonnegative function.
Then $(u_n)_{n\in\mathbb{N}}$ and $(u_n^{(\gamma^{*}+p-1)/p})_{n\in\mathbb{N}}$ are bounded sequences in $W_{\mathrm{loc}}^{1,p}(\Omega)$ and $W_0^{1,p}(\Omega)$, respectively. 
\end{enumerate}
\end{Lemma}

\begin{proof}
We prove the result only for $1<p<N$, since the other cases are analogous.

{\bf (a):} Choosing $\phi=u_n$ as a test function in the weak formulation of \eqref{approx}, we get
\begin{equation}\label{apI}
\begin{split}
\int_{\Omega}\mathcal{B}(x,\nabla u_n)\cdot\nabla u_n(x)&\,dx+\int_{\mathbb{R}^N}\int_{\mathbb{R}^N}|u_n(x)-u_n(y)|^p\,d\mu\\
&=\int_{\Omega}f_n(x)\Big(u_n(x)+\frac{1}{n}\Big)^{-\gamma(x)}u_n(x)\,dx.
\end{split}
\end{equation}
Using {\rm(H1)} in the above estimate, we obtain
\begin{equation}\label{api}
\begin{split}
&C_1\int_{\Omega}|\nabla u_n(x)|^p\,dx\leq\int_{\Omega}f_n(x)\Big(u_n(x)+\frac{1}{n}\Big)^{-\gamma(x)}u_n(x)\,dx\\
&\qquad\leq\int_{\overline{\Omega_\delta}\cap\{0<u_n\leq 1\}}f(x)u_n(x)^{1-\gamma(x)}\,dx+\int_{\overline{\Omega_\delta}\cap\{u_n>1\}}f(x)u_n(x)^{1-\gamma(x)}\,dx\\
&\qquad\qquad+\int_{\Omega\setminus\overline{\Omega_\delta}}f(x)\|c(\Omega\setminus\overline{\Omega_\delta})^{-\gamma(x)}\|_{L^\infty(\Omega)}u_n(x)\,dx\\
&\qquad\leq\|f\|_{L^1(\Omega)}+(1+\|c(\Omega\setminus\overline{\Omega_\delta})^{-\gamma(x)}\|_{L^\infty(\Omega)})\int_{\Omega}f(x)u_n(x)\,dx,
\end{split}
\end{equation}
where we have used the fact that $0<\gamma(x)\leq 1$ for every $x\in\Omega_\delta$ and the property $u_n\geq c(\Omega\setminus\overline{\Omega_\delta})>0$ in $\Omega\setminus\overline{\Omega_\delta}$ from Lemma \ref{exisapprox}. Since $f\in L^m(\Omega)$ for $m=(p^{*})'$, using H\"older's inequality and Lemma \ref{emb} in \eqref{api}, we get
\begin{equation*}
\begin{split}
C_1\|u_n\|^p&\leq \|f\|_{L^1(\Omega)}+(1+\|c(\Omega\setminus\overline{\Omega_\delta})^{-\gamma(x)}\|_{L^\infty(\Omega)})\|f\|_{L^m(\Omega)}\|u_n\|_{L^{p^{*}}(\Omega)}\\
&\leq \|f\|_{L^1(\Omega)}+C\|f\|_{L^m(\Omega)}\|u_n\|,
\end{split}
\end{equation*}
for some positive constant $C$, independent of $n$. This shows that $(u_n)_{n\in\mathbb{N}}$ is a bounded sequence in $W_0^{1,p}(\Omega)$.


{\bf (b):} Choosing $\phi=u_n$ as a test function in \eqref{approx}, proceeding as in $(a)$ above, we obtain
\begin{equation}\label{capI1}
C_1\|u_n\|^p\leq\int_{\Omega}f_n\Big(u_n+\frac{1}{n}\Big)^{-\gamma}u_n\,dx.
\end{equation}
Since $f\in L^m(\Om)$ and $(1-\gamma)m'=p^*$, using Lemma \ref{emb} in \eqref{capI1}, we get
\begin{align*}
\|u_n\|^p&\leq \int_{\Om}f(x)u_n(x)^{1-\gamma}\,dx
\leq \|f\|_{L^m(\Om)}\left(\int_{\Om}u_n(x)^{(1-\gamma)m'}\right)^\frac{1}{m'}\\
&=\|f\|_{L^m(\Om)}\left(\int_{\Om}u_n(x)^{p^*}\,dx\right)^\frac{1-\gamma}{p^*}
\leq C\|f\|_{L^m(\Om)}\|u_n\|^{1-\gamma},
\end{align*}
for some constant $C$, independent of $n$. Therefore, we have
$\|u_n\|\leq C$ for every $n\in\mathbb N$.
Thus $(u_n)_{n\in\mathbb{N}}$ is a bounded sequence in $W_0^{1,p}(\Om)$.

{\bf (c):} Choosing $\phi=u_n$ as a test function in \eqref{approx}, proceeding as in $(a)$ above, we obtain
\begin{equation}\label{capI}
\begin{split}
C_1\|u_n\|^p\leq\|f\|_{L^1(\Om)}.
\end{split}
\end{equation}
This shows that $(u_n)_{n\in\mathbb{N}}$ is a bounded sequence in $W_0^{1,p}(\Om)$.

{\bf (d):} Note that by Lemma \ref{exisapprox}, we know that $u_n\in W_0^{1,p}(\Omega)\cap L^\infty(\Omega)$. Therefore, choosing $u_n^{\gamma^*}$ as a test function in the weak formulation of \eqref{approx}, we obtain
\begin{equation}\label{apiGrt1}
\begin{split}
\int_{\Omega}\mathcal{B}(x,\nabla u_n)\cdot\nabla (u_n(x)^{\gamma^*})&\,dx
+\int_{\mathbb{R}^N}\int_{\mathbb{R}^N}
\mathcal{A}u_n(x,y) (u_n(x)^{\gamma^*}-u_n(y)^{\gamma^*})\,d\mu\\
&=\int_{\Omega}f_n(x)\Big(u_n(x)+\frac{1}{n}\Big)^{-\gamma(x)}u_n(x)^{\gamma^*}\,dx.
\end{split}
\end{equation}
We observe that by Lemma \ref{BrPrapp} the nonlocal integral in the above estimate \eqref{apiGrt1} becomes nonnegative. Thus applying {\rm(H1)} in \eqref{apiGrt1} gives
\begin{equation}\label{apigrt1}
\begin{split}
&C_1\gamma^*\Big(\frac{p}{\gamma^{*}+p-1}\Big)^p\int_{\Omega}\Big|\nabla\big(u_n(x)^\frac{\gamma^{*}+p-1}{p}\big)\Big|^{p}\,dx\\
&\qquad\leq\int_{\Omega}f_n(x)\Big(u_n(x)+\frac{1}{n}\Big)^{-\gamma(x)}u_n(x)^{\gamma^*}\,dx\\
&\qquad\leq\int_{\overline{\Omega_\delta}\cap\{0<u_n\leq 1\}}f(x)u_n(x)^{\gamma^{*}-\gamma(x)}\,dx
+\int_{\overline{\Omega_\delta}\cap\{u_n>1\}}f(x)u_n(x)^{\gamma^{*}-\gamma(x)}\,dx\\
&\qquad\qquad+\int_{\Omega\setminus\overline{\Omega_\delta}}f(x)\|c(\Omega\setminus\overline{\Omega_\delta})^{-\gamma(x)}\|_{L^\infty(\Omega)}u_n(x)^{\gamma^*}\,dx\\
&\qquad\leq\|f\|_{L^1(\Omega)}+(1+\|c(\Omega\setminus\overline{\Omega_\delta})^{-\gamma(x)}\|_{L^\infty(\Omega)})\int_{\Omega}f(x)u_n(x)^{\gamma^*}\,dx\\
&\qquad\leq \|f\|_{L^1(\Omega)}+(1+\|c(\Omega\setminus\overline{\Omega_\delta})^{-\gamma(x)}\|_{L^\infty(\Omega)})\int_{\Omega}f(x)\Big(u_n(x)^{\frac{\gamma^{*}+p-1}{p}}\Big)^\frac{p\gamma^*}{\gamma^*+p-1}\,dx\\
&\qquad\leq\|f\|_{L^1(\Omega)}+c
\|f\|_{L^m(\Omega)}\Big\|u_n^{\frac{\gamma^{*}+p-1}{p}}\Big\|^\frac{p^*}{m'},
\end{split}
\end{equation}
for some positive constant $c$. To obtain the above estimate, we have used Lemma \ref{emb} along with the hypothesis $\|\gamma\|_{L^\infty(\Omega_\delta)}\leq\gamma^*$, the fact $u_n\geq c(\Omega')>0$ from Lemma \ref{exisapprox} and also the property $u_{n}^{(\gamma^*+p-1)/p}\in W_0^{1,p}(\Omega)$, which is true since $u_n\in W_0^{1,p}(\Omega)\cap L^\infty(\Omega)$. Since $p>\frac{p^*}{m'}$, it follows from \eqref{apigrt1} that $(u_{n}^{(\gamma^*+p-1)/p})_{n\in\mathbb{N}}$ is uniformly bounded in $W_0^{1,p}(\Omega)$. Using this fact and H\"older's inequality, for any $\Omega'\Subset\Omega$, we observe that 
\begin{equation}\label{lpbd}
\int_{\Omega'}u_{n}(x)^p\,dx
\leq\left(\int_{\Omega'}u_{n}(x)^{\gamma^{*}+p-1}\,dx\right)^\frac{p}{\gamma^{*}+p-1}|\Om|^\frac{\gamma^{*}-1}{\gamma^{*}+p-1}\leq C, 
\end{equation}
for some positive constant $C$ independent of $n$ and
\begin{equation}\label{gradbd}
\begin{split}
\int_{\Omega'}\Big|\nabla\big(u_{n}(x)^\frac{\gamma^{*}+p-1}{p}\big)\Big|^p\,dx
&=\Big(\frac{\gamma^{*}+p-1}{p}\Big)^p\int_{\Omega'}u_{n}(x)^{\gamma^{*}-1}|\nabla u_n(x)|^p\,dx\\
&\geq c(\Omega')^{\gamma^*-1}\Big(\frac{\gamma^{*}+p-1}{p}\Big)^p\int_{\Omega'}|\nabla u_n(x)|^p\,dx.
\end{split}
\end{equation}
In the last line above, we have again used the fact that $u_n\geq c(\Omega')>0$ for every $\Omega'\Subset\Omega$ for some positive constant $c(\Omega')$, independent of $n$. Therefore, the estimates \eqref{lpbd} and \eqref{gradbd} yields the uniform boundedness of $(u_n)_{n\in\mathbb{N}}$ in $W^{1,p}_{\mathrm{loc}}(\Omega)$. This completes the proof.

{\bf (e):} Let $\gamma(x)=\gamma>1$. Again, note that by Lemma \ref{approx}, we know that $u_n\in W_0^{1,p}(\Omega)\cap L^\infty(\Omega)$ and as in $(c)$ above, choosing $u_n^{\gamma}$ as a test function in \eqref{approx}, we obtain
\begin{equation}\label{apigrt}
\begin{split}
&C_1\gamma\Big(\frac{p}{\gamma+p-1}\Big)^p\int_{\Omega}\Big|\nabla \big(u_n(x)^\frac{\gamma+p-1}{p}\big)\Big|^{p}\,dx\\
&\qquad=\int_{\Omega}f_n(x)\Big(u_n(x)+\frac{1}{n}\Big)^{-\gamma}u_n(x)^{\gamma}\,dx\leq\|f\|_{L^1(\Om)}.
\end{split}
\end{equation}
Thus, proceeding similarly as in (c) above, the result follows.
\end{proof}

Next, we establish the pointwise convergence of the gradient of the approximate solutions $(u_n)_{n\in\mathbb{N}}$ found in Lemma \ref{approx}. Here we adapt the proof of \cite[Theorem 7.1]{GU} to the variable exponent and anisotropic case.
\begin{Lemma}\label{grad}
Let $1<p<\infty$. Suppose that $(u_n)_{n\in\mathbb{N}}$ is the sequence of approximate solutions for the problem \eqref{approx} given by Lemma \ref{exisapprox} and $u$ is the pointwise limit of $(u_n)_{n\in\mathbb{N}}$. 
Let $\gamma$ and $f$ satisfy any one of the conditions $(a)-(e)$ in Lemma \ref{apun}. 
Then there exists a subsequence, denoted again by $(u_n)_{n\in\mathbb{N}}$, such that $\nabla u_n\to \nabla u$ pointwise almost everywhere in $\Om$ as $n\to\infty$.
\end{Lemma}
\begin{proof}
Lemma \ref{apun} implies that 
the sequences $(u_n)_{n\in\mathbb{N}}$ and $(u_n^{(\gamma^{*}+p-1)/p})_{n\in\mathbb{N}}$ are uniformly bounded in $W^{1,p}_{\mathrm{loc}}(\Om)$ and $W_0^{1,p}(\Om)$, respectively.
Therefore, we have $u_n\rightharpoonup u$ weakly in $W^{1,p}_{\mathrm{loc}}(\Om)$ and 
$u_n\to u$ in $L^p_{\mathrm{loc}}(\Om)$ as $n\to\infty$.
Let $K\subset\Omega$ be a compact set and consider a function $\phi_K\in C_c^{1}(\Omega)$ such that $\mathrm{supp}\,\phi_K=\Omega'$, $0\leq\phi_K\leq 1$ in $\Omega$ and $\phi_K\equiv 1$ in $K$. For $\mu>0$, we define the truncation $T_{\mu}:\mathbb{R}\to\mathbb{R}$ by
\begin{equation}\label{trun}
T_{\mu}(s)=
\begin{cases}
s,\quad\text{ if }|s|\leq\mu,\\
\mu\frac{s}{|s|},\quad\text{ if }|s|>\mu.
\end{cases}
\end{equation}
Choosing $v_n=\phi_K T_{\mu}((u_n-u))\in W_0^{1,p}(\Omega)$ as a test function in \eqref{approx}, we obtain
\begin{equation}\label{gs1eqn1}
\begin{split}
I+J=R,
\end{split}
\end{equation}
where
$$
I=\int_{\Omega}\mathcal{B}(x,\nabla u_n)\cdot\nabla v_n(x)\,dx,
$$
$$
J=\int\limits_{\mathbb{R}^N}\int\limits_{\mathbb{R}^N}\mathcal{A}u_n(x,y)\big(v_n(x)-v_n(y)\big)\,d\mu
\quad\text{and}
\quad R=\int_{\Omega}\frac{f_n(x)}{\big(u_n(x)+\frac{1}{n}\big)^{\gamma(x)}} v_n(x)\,dx.
$$

\textbf{Estimate for $I$:} 
We observe that
\begin{equation}\label{estI}
\begin{split}
I&=\int_{\Omega}\mathcal{B}(x,\nabla u_n)\cdot\nabla v_n(x)\,dx\\
&=\int_{\Omega}\phi_K(x)\big(\mathcal{B}(x,\nabla u_n)-\mathcal{B}(x,\nabla u)\big)\cdot\nabla T_{\mu}(u_n(x)-u(x))\,dx\\
&\qquad+\int_{\Omega}\phi_K(x)\mathcal{B}(x,\nabla u)\cdot\nabla T_{\mu}(u_n(x)-u(x))\,dx\\
&\qquad+\int_{\Omega}T_{\mu}(u_n(x)-u(x))\mathcal{B}(x,\nabla u_n)\cdot\nabla\phi_K(x)\,dx\\
&=I_1+I_2+I_3.
\end{split}
\end{equation}

\textbf{Estimate for $I_2$:} Since $u_n\rightharpoonup u$ weakly in $W^{1,p}_{\mathrm{loc}}(\Om)$ as $n\to\infty$, we have $T_{\mu}(u_n-u)\rightharpoonup 0$ weakly in $W^{1,p}_{\mathrm{loc}}(\Omega)$ as $n\to\infty$. By {\rm(H1)}, we get
\begin{equation}\label{estI_2}
\lim_{n\to\infty}I_2
=\lim_{n\to\infty}\int_{\Omega}\phi_K(x)\mathcal{B}(x,\nabla u)\cdot\nabla T_{\mu}(u_n(x)-u(x))\,dx=0.
\end{equation}

\textbf{Estimate for $I_3$:} For $\Omega'=\mathrm{supp}\,\phi_K$, by H\"older's inequality, using the uniform boundedness of $(u_n)_{n\in\mathbb{N}}$ in $W^{1,p}(\Omega')$ and the hypothesis {\rm(H1)}, we have
\begin{equation}\label{estI_3}
\begin{split}
|I_3|
&=\left|\int_{\Omega}T_{\mu}(u_n(x)-u(x))\mathcal{B}(x,\nabla u_n)\cdot\nabla\phi_K(x)\,dx\right|\\
&\leq C_2\int_{\Omega}|\nabla\phi_K(x)|\,|u_n(x)-u(x)|\,|\nabla u_n(x)|^{p-1}\,dx\\
&\leq C_2\Vert\nabla\phi_K\Vert_{L^\infty(\Omega')}\left(\int_{\Omega'}|\nabla u_n(x)|^p\,dx\right)^\frac{p-1}{p}\Vert u_n-u\Vert_{L^p(\Omega')}\\
&\leq C\Vert u_n-u\Vert_{L^p(\Omega')},
\end{split}
\end{equation}
for some constant $C$, independent of $n$. Since $u_n\to u$ in $L^p_{\mathrm{loc}}(\Om)$, by \eqref{estI_3} as $n\to\infty$, we obtain
\begin{equation}\label{estI_3final}
\lim_{n\to\infty}I_3
=\lim_{n\to\infty}\int_{\Omega}T_{\mu}(u_n(x)-u(x))\mathcal{B}(x,\nabla u_n)\cdot\nabla\phi_K(x)\,dx=0.
\end{equation}
Therefore, by \eqref{estI_2} and \eqref{estI_3final} in \eqref{estI}, we arrive at
\begin{equation}\label{estIfinal}
\begin{split}
\limsup_{n\to\infty}I
&=\limsup_{n\to\infty}I_1\\
&=\limsup_{n\to\infty}\int_{\Omega}\phi_K(x)\big(\mathcal{B}(x,\nabla u_n)-\mathcal{B}(x,\nabla u)\big)\cdot\nabla T_{\mu}(u_n(x)-u(x))\,dx.
\end{split}
\end{equation}

\textbf{Estimate for $J$:}  Since $u_n\rightharpoonup u$ weakly in $W^{1,p}_{\mathrm{loc}}(\Om)$ as $n\to\infty$, we may proceed along the lines of the proof of the estimate $(A.36)$ in \cite[pages 29-33]{GU}, it follows that
\begin{equation}\label{estJfinal}
\lim_{n\to\infty}J\geq 0.
\end{equation}

\textbf{Estimate for $R$:} Recalling that $\mathrm{supp}\,\phi_K=\Omega'$, by Lemma \ref{approx}, there exists a constant $C(\Omega')>0$, independent of $n$, such that $u\geq C(\Omega')>0$ in $\Om$. Hence, we have
\begin{equation}\label{estK}
R=\int_{\Om}\frac{f_n(x)}{\big(u_n(x)+\frac{1}{n}\big)^{\gamma(x)}}v_n(x)\,dx
\leq\frac{\|f\|_{L^1(\Om)}}{\|C(\Omega')^{-\gamma}\|_{L^\infty(\Omega')}}\mu.
\end{equation}
Therefore, for every fixed $\mu>0$, using \eqref{estIfinal}, \eqref{estJfinal} and \eqref{estK} in \eqref{gs1eqn1}, we obtain
\begin{equation}\label{gest}
\begin{split}
\limsup_{n\to\infty}\int_{K}\big(\mathcal{B}(x,\nabla u_n)-\mathcal{B}(x,\nabla u)\big)\cdot\nabla T_{\mu}(u_n(x)-u(x))\,dx\leq C\mu,
\end{split}
\end{equation}
for some constant $C=C\big(\Omega',\|f\|_{L^1(\Om)}\big)>0$. Taking into account \eqref{gest} along with the hypotheses on $\mathcal{B}(x,\nabla u)$, 
the rest of the proof follows from Step 2 in \cite[page 586]{BocMur}.
\end{proof}

\subsection{Proof of the existence results}

\begin{proof}[Proof of Theorem \ref{thm2}]
{\bf (a):} Let $1<p<N$ and $f\in L^m(\Om)\setminus\{0\}$ be a nonnegative function with $m=(p^{*})'$. 
By Lemma \ref{exisapprox}, for every $n\in\mathbb{N}$, there exists $u_n\in W_0^{1,p}(\Omega)$ such that 
\begin{equation}\label{leqn1}
\begin{split}
&\int_{\Om}\mathcal{B}(x,\nabla u_n)\cdot\nabla\phi(x)\,dx
+\int_{\mathbb{R}^N}\int_{\mathbb{R}^N}\mathcal{A}u_n(x,y)(\phi(x)-\phi(y))\,d\mu\\
&\qquad=\int_{\Om}f_n(x)\Big(u_n^{+}(x) +\frac{1}{n}\Big)^{-\gamma(x)}\,\phi(x)\,dx
\end{split}
\end{equation}
for every $\phi\in C_c^{1}(\Om)$.
By Lemma \ref{apun}, we have $\|u_n\|\leq C$, $n\in\mathbb N$, for some positive constant $C$, independent of $n$ and $u_{n+1}\geq u_n$ for every $n\in\mathbb{N}$ from Lemma \ref{exisapprox}. Therefore, there exists $u\in W_0^{1,p}(\Omega)$ such that $0<u_n\leq u$ almost everywhere in $\Om$ and up to a subequence $u_n\to u$ pointwise almost everywhere in $\Omega$ and weakly in $W_0^{1,p}(\Omega)$. Moreover, by Lemma \ref{exisapprox}, for every $\Omega'\Subset\Omega$, there exists a positive constant $c(\Omega')$, independent of $n$ such that $u\geq u_n\geq c(\Omega')>0$ in $\Omega'$, which gives
\begin{equation}\label{rhs}
\left|f_n(x)\Big(u_n(x)+\frac{1}{n}\Big)^{-\gamma(x)}\phi\right|
\leq \|c_{\mathrm{supp}\,\phi}^{-\gamma(x)}\,\phi\|_{L^\infty(\Omega)}f(x)\in L^1(\Omega),
\end{equation}
for every $\phi\in C_c^{1}(\Omega)$.
Taking into account \eqref{rhs}, we may apply the dominated convergence theorem to obtain
\begin{equation}\label{rhslim}
\lim_{n\to\infty}\int_{\Omega}f_n(x)\Big(u_n(x)+\frac{1}{n}\Big)^{-\gamma(x)}\phi(x)\,dx
=\int_{\Omega}f(x)u(x)^{-\gamma(x)}\phi(x)\,dx
\end{equation}
for every $\phi\in C_c^{1}(\Omega)$.

Moreover, by Lemma \ref{grad}, up to a subsequence $\nabla u_n\to \nabla u$ pointwise almost everywhere in $\Omega$ as $n\to\infty$ and proceeding  as in the proof of \eqref{d1} above, we obtain
\begin{equation}\label{lhs}
\lim_{n\to\infty}\int_{\Omega}\mathcal{B}(x,\nabla u_n)\cdot\nabla\phi(x)\,dx
=\int_{\Omega}\mathcal{B}(x,\nabla u)\cdot\nabla\phi(x)\,dx
\end{equation}
for every $\phi\in C_c^{1}(\Omega)$.

Since $\phi\in C_c^1(\Om)$ and $(u_n)_{n\in\mathbb{N}}$ is uniformly bounded in $W_0^{1,p}(\Om)$, by Lemma \ref{locnon1}
$$
\frac{|u_n(x)-u_n(y)|^{p-2}\big(u_n(x)-u_n(y)\big)}{|x-y|^\frac{N+ps}{p'}},\quad n\in\mathbb N,
$$
is uniformly bounded in $L^{p'}(\mathbb{R}^N\times\mathbb{R}^N)$ and
$$
\frac{\phi(x)-\phi(y)}{|x-y|^\frac{N+ps}{p}}\in L^p(\mathbb{R}^N\times\mathbb{R}^N).
$$
Therefore, by the weak convergence, we have
\begin{equation}\label{nonloclim}
\lim_{n\to\infty}\int_{\mathbb{R}^N}\int_{\mathbb{R}^N}\mathcal{A}u_n(x,y)(\phi(x)-\phi(y))\,d\mu
=\int_{\mathbb{R}^N}\int_{\mathbb{R}^N}\mathcal{A}u(x,y)(\phi(x)-\phi(y))\,d\mu.
\end{equation}
Thus using \eqref{rhslim}, \eqref{lhs} and \eqref{nonloclim} in \eqref{leqn1}, the result follows. The proof for $p\geq N$ is analogous.

{\bf (b):} By Lemma \ref{apun} (d), the sequences $(u_n)_{n\in\mathbb{N}}$ and $(u_n^{(\gamma^{*}+p-1)/p})_{n\in\mathbb{N}}$ are uniformly bounded in $W^{1,p}_{\mathrm{loc}}(\Omega)$ and $W^{1,p}_{0}(\Omega)$, respectively. Then, analogous to the proof of (a) above, the result follows. 
\end{proof}

\begin{proof}[Proof of Theorem \ref{cthm1}]  
By [Lemma \ref{apun} $(b)$], the sequences $(u_n)_{n\in\mathbb{N}}$ is uniformly bounded in $W^{1,p}_0(\Omega)$. Then, analogous to the proof of (a) above, the result follows. 
\end{proof}

\begin{proof}[Proof of Theorem \ref{cthm2}] 
By Lemma \ref{apun} (c), the sequences $(u_n)_{n\in\mathbb{N}}$ is uniformly bounded in $W^{1,p}_0(\Omega)$. Then, analogous to the proof of (a) above, the result follows. 
\end{proof}

\begin{proof}[Proof of Theorem \ref{cthm3}] 
By Lemma \ref{apun} (e), the sequences $(u_n)_{n\in\mathbb{N}}$ and $(u_n^{(\gamma^{*}+p-1)/p})_{n\in\mathbb{N}}$ are bounded in $W^{1,p}_{\mathrm{loc}}(\Omega)$ and $W^{1,p}_{0}(\Omega)$, respectively. Then, analogous to the proof of (a) above, the result follows. 
\end{proof}

\noindent {\textsf{Prashanta Garain\\
Department of Mathematical Sciences,\\
Indian Institute of Science Education and Research Berhampur,\\
Berhampur, Odisha 760010, India},\\
\textsf{e-mail}: pgarain92@gmail.com\\

\noindent {\textsf{Wontae Kim\\
Department of Mathematics, Aalto University,\\
P.O. BOX 11100, 00076, Aalto, Finland},\\
\textsf{e-mail}: wontae.kim@aalto.fi\\

\noindent {\textsf{Juha Kinnunen\\
Department of Mathematics, Aalto University,\\
P.O. BOX 11100, 00076, Aalto, Finland},\\
\textsf{e-mail}: juha.k.kinnunen@aalto.fi\\

\end{document}